\documentclass[12pt]{amsart}

\newtheorem{thm}{Theorem}

\newtheorem{prop}{Proposition}

\DeclareMathOperator{\der}{der}

\DeclareMathOperator{\Lie}{Lie}

\DeclareMathOperator{\tr}{tr}

\DeclareMathOperator{\Adj}{Ad}

\DeclareMathOperator{\reg}{reg}

\DeclareMathOperator{\Span}{span}

\newcommand{\gm}{\gamma}

\newcommand{\half}{\frac{1}{2}}

\newcommand{\BB}{\mathcal B}  
\newcommand{\CC}{\ensuremath{{\mathbb C}}}
\newcommand{\Ch}{\ensuremath{{\mathcal C}}}

\newcommand{\II}{\ensuremath{{\mathfrak a}}}

\newcommand{\PP}{\ensuremath{{\mathcal P}}}

\newcommand{\R}{\ensuremath{{\mathbb R}}}

\newcommand{\ZZ}{\ensuremath{{\mathbb Z}}}

\renewcommand{\lim}[1]{\underset{#1}{\operatorname{lim}}}

\begin{document}

\title{On Arthur's $\Phi$-Function}
\author{Steven Spallone}
\address{Max-Planck-Institut f\"{u}r Mathematik, Bonn, Germany}
\date{\today}
\email{spallone@mpig-bonn.mpg.de}


\begin{abstract}

Write $\Theta^E$ for the stable character associated to a finite 
dimensional representation $E$ of a connected real reductive group $G$.  
Let $M$ be the centralizer of a maximal torus $T$, and denote by 
$\Phi_M(\gm,\Theta^E)$ Arthur's extension of $ |D_M^G(\gm)|^{\half} 
\Theta^E(\gm)$ to $T(\R)$.  In this paper we give a simple 
explicit expression for 
$\Phi_M(\gm,\Theta^E)$, when $\gm$ is elliptic in $G$.
\end{abstract}

\maketitle

\section{Introduction}
Let $G$ be a connected reductive group over $\R$, and $T$ a 
maximal torus in $G$.  Assume that $G$ has a discrete series of 
representations.  Let $A$ 
be the split part of $T$, and $M$ the centralizer of $A$ in $G$.  It is a Levi 
subgroup of $G$ 
containing $T$.  Let $E$ be a finite-dimensional representation of 
$G(\CC)$, 
and consider the packet $\Pi_E$ of discrete series representations $\pi$ of 
$G(\R)$ which have the same infinitesimal and central characters as $E$.  
Write $\Theta_{\pi}$ for the character of $\pi$, and put
\[ \Theta^E=(-1)^{q(G)} \sum_{\pi \in \Pi_E} \Theta_{\pi}. \]

\noindent Here $q(G)$ is half the dimension of the symmetric space 
associated with $G$.
Note that $\Theta^E(\gm)$ will not extend to all elements $\gm \in T(\R)$, 
in particular to $\gm=1$. 
Define the number $D_M^G(\gm)$ by
\[ D_M^G(\gm) = \det(1-\Adj(\gm); \Lie(G)/ \Lie(M)). \]
Then a result of Arthur and Shelstad [1] states that the 
function
\[ \gm \mapsto |D_M^G(\gm)|^{\half} \Theta^E(\gm), \]

\noindent defined on the set of regular elements $T_{\reg}(\R)$
extends continuously to $T(\R)$.  We denote this extension by 
$\Phi_M(\gm,\Theta^E)$.  These quantities give the contribution 
from the real place to the $L^2$-Lefschetz numbers of Hecke operators in 
[1] and [2].
An expression for $\Phi_M(\gm,\Theta^E)$ as essentially a sum over elements 
in the Weyl group $W$ of $T$ in $G$ appears in the proof of Lemma 4.1 in 
[2].  Although this expression suffices to prove the lemma, it can be 
considerably refined when $\gm$ is in the maximal elliptic subtorus 
$T_e(\R)$ of $T(\R)$.

\bigskip

The following theorem is proved in section 4.

\bigskip

{\bf Theorem 1.} {\it If $\gm \in T_e(\R)$, then \it}
\[ \Phi_M(\gm,\Theta^E)=(-1)^{q(L)}\cdot |W_L| \cdot  
\sum_{\omega \in W^{LM}}
\varepsilon(\omega) \cdot
\tr(\gm;V^M_{\omega(\lambda_B+\rho_B)-\rho_B}). \]

Here we write $L$ for the centralizer of $T_c$ in $G$, where $T_c$ 
is the maximal compact subtorus of $T$.  Also write $W_L$ 
and $W_M$ for the Weyl groups of $T$ in $L$ and $M$.  The latter
are subgroups of $W$ which commute and have trivial intersection.
Here $W^{LM}$ is a certain set of representatives for the cosets  
$(W_L \times W_M) \backslash W$.  It is defined explicitly in section 5. 
We write $\varepsilon$ 
for the sign character of $W$.  Finally by 
$V^M_{\omega(\lambda_B+\rho_B)-\rho_B}$ we
denote the irreducible finite-dimensional representation of $M$,
with highest weight $\omega(\lambda_B+\rho_B)-\rho_B$, where $\lambda_B$ 
is the $B$-dominant highest weight of $E$.

\bigskip

In particular, we obtain the extremely simple expression,
\[ \Phi_A(1,\Theta^E)=(-1)^{q(G)} \cdot |W|, \]

\noindent in the case of a split torus $T=A$.

\bigskip

We now describe the organization of this paper.

\bigskip

In section 2, we spell out the relationship between the root 
systems of $G$, $L$, and $M$.  There are two distinct systems of chambers in 
$X_*(A) \otimes_{\ZZ} \R$ obtained from these root systems which are important to 
understand.

\bigskip

In section 3, we take the aforementioned lemma a step further to 
express $\Phi_M(\gm,\Theta^E)$ explicitly as a linear combination of 
characters.  (Actually we do the computation for any stable virtual 
character $\Theta$, as it is no more difficult.)  The sum over $W$ simplifies 
to a sum over Kostant representatives $W^M$. 

\bigskip

In section 4, in which we deal specifically with 
$\Phi_M(\gm,\Theta^E)$, we 
distill out the action of $W_L$.  A sum over $W^{LM}$ remains.  At a key step we use 
a result of section 5, the computation of an alternating sum of stable 
discrete 
series constants.

\bigskip

In section 5, we prove the result mentioned above, in the 
context of abstract root systems.  It is independent from the rest of the paper.

\bigskip

I am indebted to my advisor Robert Kottwitz for suggesting the problem and 
many useful comments.  I also thank Christian Kaiser for some helpful conversations.  

This project was carried out during a stay at the Max-Planck-Institut 
f\"{u}r Mathematik in Bonn, and I am grateful to the Institut for its 
support and hospitality.

\section{$L$-chambers and $\PP$-chambers}

Let $G$ be a connected reductive group over $\R$, and $T$ a 
maximal torus of $G$.  
Assume that $G$ has a discrete series, or equivalently, that $G$ has an 
elliptic maximal torus.
\bigskip

Write $T_c$, respectively $A$, for the maximal compact, 
resp. split, 
subtori of $T$ with centralizers $L$, resp. $M$, in $G$. 
Write $R$ for the root system of $T$ in $G$, and $R_L$, 
resp. $R_M$, for the set of roots of $T$ in $L$, resp. $M$.  Then $R_L$ is the 
subset of $R$ consisting of real roots, and $R_M$ is the subset of imaginary roots. 
Write $W_L$ and $W_M$ for the respective Weyl groups.  They are commuting subgroups 
of $W$ with trivial intersection.  Note that $W_L$ fixes each root in $R_M$.
\bigskip

$A$ is contained as a split maximal torus in $L_{\der}$, the 
derived group of $L$, and we 
may identify $R_L$ with the set of roots of $A$ in $L_{\der}$.
\bigskip

Write $\II_M$ for $X_*(A) \otimes_{\ZZ} \R$. 
For any $\alpha \in R \backslash R_M$ the root hyperplane 
$H_{\alpha}$ of 
$X^*(T)_{\R}:=X_*(T) \otimes \R$ gives a hyperplane in $\II_M$.  Let us 
consider two kinds of chambers in $\II_M$ obtained from these.
Define $\PP$-chambers to be those obtained by deleting from $\II_M$ all 
the hyperplanes $H_{\alpha}$, with $\alpha \in R \backslash R_M$. 
Define $L$-chambers to be those obtained by deleting all the $H_{\alpha}$ 
with $\alpha \in R_L$.  The latter are the Weyl chambers for $A$ in 
$L_{\der}$; therefore $W_L$ acts simply transitively on them.
\bigskip

Observe that $R_L \subset (R \backslash R_M)$.  Any additional hyperplanes 
coming from roots in $R \backslash (R_L \cup R_M)$ divide the $L$-chambers
into $\PP$-chambers.  Thus every $\PP$-chamber is contained in a unique $L$-chamber.
\bigskip

Write $\PP(M)$ for the set of parabolic subgroups of $G$ 
admitting $M$ as a Levi 
component.  There is a one-to-one correspondence between $\PP(M)$ and the set of 
$\PP$-chambers in $\II_M$, obtained as follows:  for $P=MN \in \PP(M)$, the 
corresponding $\PP$-chamber is 
\[ \II_P^+= \{x \in \II_M : \langle \alpha,x \rangle > 0 \text{, for all } \alpha \in 
R_N \}, \]

\noindent where $R_N$ denotes the set of roots of $T$ in $\Lie(N)$.
\bigskip

Recall that the set of $L$-chambers is in bijection with the set 
of Borel 
subgroups of $L$ containing $T$, or equivalently the set of positive root systems
$R_L^+$ in the root system $R_L$.
\bigskip

Now let $C_P$ be a $\PP$-chamber, and let $P=MN$ be the 
corresponding element of 
$\PP(M)$.  It is easy to see that $R_N \cap R_L$ is a positive system in 
$R_L$, 
and this corresponds to an $L$-chamber $C_L$.  Thus we have defined a map $C_P 
\mapsto C_L$ from the set of $\PP$-chambers to the set of $L$-chambers.  It is 
the obvious one which associates to $C_P$ the unique $L$-chamber 
containing $C_P$.

\section{A Linear Combination of Characters}

A stable virtual character is a finite $\ZZ$-linear combination 
$\Theta$ of 
characters $\Theta_{\pi}$ so that
\[ \Theta(\gm)=\Theta(\gm^{\prime}) \]
\noindent whenever $\gm$ and $\gm^{\prime}$ are regular, stably conjugate 
elements of 
$G(\R)$.

In Lemma 4.1 of [2], it is proved that for a stable virtual 
character $\Theta$ on 
$G(\R)$, the function
\[ \gm \mapsto |D^G_M(\gm)|^{\half}\Theta(\gm) \]
\noindent on $T_{\reg}(\R)$ extends continuously to $T(\R)$.  A key 
ingredient of the proof is the 
fact that the expression at the bottom of page 497 is a linear combination 
of irreducible finite-dimensional representations of $M$.  
In this section we will compute explicitly the 
coefficients and the representations involved, in the case where the 
element $a$ appearing in the proof is equal to $1$.

\bigskip
We translate the set-up of the proof in [2] as follows.  
We take $\Gamma$ to be the identity component of $T(\R)$.  The root 
system $R_{\Gamma}$ is then simply $R_L$.  Fix a 
positive root system $R_L^+$ in $R_L$, and let $C$ be the corresponding 
$L$-chamber in $\II_M$.  
We then choose a parabolic subgroup $P=MN$ so that $R_L \cap R_N \subseteq R_L^+$.
Note that $R_L \cap R_N$ is also a system of positive roots, so this condition is 
equivalent to having $R_L \cap R_N=R_L^+$.  Thus we simply require that 
the $\PP$-chamber corresponding to $P$ be contained in $C$.
\bigskip

Although at the end of our computations we will allow $\gm$ to 
be nonregular, 
we choose now $\gm$ to be a regular element of $\Gamma= T_c(\R) \cdot 
\exp(\bar{C})$.
\bigskip

The expression is
\begin{equation}
 \sum_B m(B) \frac{\Delta_P(\gm) \cdot \lambda_B(\gm)}{\Delta_B(\gm)}. 
\end{equation}

\noindent The sum ranges over Borels containing $T$, which correspond to 
elements of 
$W$.

\bigskip

\noindent Here $\lambda_B$ is the $B$-dominant highest weight of $E$,  
\[ \Delta_B=\prod_{\alpha >0} (1-\alpha^{-1}) \text{, and } 
\Delta_P=\prod_{\alpha \in R_N}(1-\alpha^{-1}). \]

Fix now a Borel $B$ of $G$ with $T \subseteq B \subseteq P$, for 
the rest of 
this paper.

\bigskip

Recall the set of Kostant representatives $W^M$ 
for the Weyl group $W_M$ of $M$, relative to B.
It is the set $\{ w \in W| w^{-1}R_M^+ \subset R^+ \}$.

\bigskip

If $w \in W$, write $w * B$ for $wBw^{-1}$.

\bigskip

We will use the observation that for $\omega \in W^M, (\omega * 
B)_M=B_M$.  
Indeed, if $\alpha \in R^+ \cap R_M$, then $\omega^{-1} \alpha \in R^+$, which implies that $\alpha \in   
\omega R^+ \cap R_M$.

\bigskip

Our sum (1) breaks up as

\begin{equation}
\label{break}
\sum_{\omega \in W^M} m(\omega * B) \cdot \Delta_P(\gm) \cdot \sum_{w_M \in W_M} 
\frac{w_M (\omega\lambda_B)(\gm)}{\Delta_{w_M \omega * B}(\gm)}.
\end{equation}
We would prefer the denominator inside the sum to be $\Delta_{
w_M * B_M}(\gm)$.
(Recall that $B_M=B \cap M$.)  Note that $\Delta_P \cdot \Delta_{B_M} = \Delta_B$, since
$R^+$ is the disjoint union of $R_M^+$ and $R_N$.

\bigskip

So we consider the quantity

\begin{equation}
\label{quantity}
 \frac{ \Delta_P \cdot \Delta_{w_M * B_M}}{\Delta_{w_M \omega * B}}=
\frac{ \Delta_B \cdot \Delta_{w_M * B_M}}{\Delta_{B_M} \cdot \Delta_{w_M \omega * 
B}} . 
\end{equation}

\noindent Observe that if $\BB$ is a Borel, $\Delta_{\BB}=\delta_{\BB} 
\cdot 
\rho_{\BB}^{-1}$, 
where $\delta_{\BB}=\prod_{\alpha > 0} 
(\alpha^{\half}-\alpha^{-\half})$ and $\rho_{\BB}$ is the usual half sum of positive 
roots.  
Since $\delta_{w * \BB}=\varepsilon(w)\delta_{\BB}$, we compute that 
\[ \frac{\Delta_{w * \BB}}{\Delta_{\BB}}=\varepsilon(w) \cdot 
(\rho_{\BB}-w\rho_{\BB}). \]

Thus \eqref{quantity} becomes
\[ \varepsilon(\omega) (w_M(\omega \rho_B-\rho_{B_M})-\rho_B+\rho_{B_M}). \]

Next observe that for $w_M \in W_M$,
\[ w_M(\rho_B-\rho_{B_M})=\rho_B-\rho_{B_M}. \]

Indeed, the roots of $R^+$ not in $R_M^+$ are in $R_N$, and are 
thus normalized by 
$W_M$.  So the above expression simplifies to
\[ \varepsilon(\omega) \cdot w_M(\omega \rho_B-\rho_B). \]

We can therefore rewrite \eqref{break} as

\begin{equation}
\label{rewrite}
\sum_{\omega \in W^M} m(\omega * B) \cdot \varepsilon(\omega) \cdot \sum_{w_M \in W_M}
\frac{w_M (\omega(\lambda_B+\rho_B)-\rho_B)(\gm)}{\Delta_{w_M * B_M}(\gm)}.
\end{equation}

Since $\omega$ is a Kostant representative, the 
weight 
$\omega(\lambda_B+\rho_B)-\rho_B$ is positive for $B_M$, 
and we may use the Weyl character formula to rewrite this as

\begin{equation}
\label{Weyl}
\sum_{\omega \in W^M} m(\omega * B) \cdot \varepsilon(\omega) \cdot \tr(\gm; V^M_{\omega(\lambda_B+\rho_B)-\rho_B}).
\end{equation}

\noindent Here $V^M_{\omega(\lambda_B+\rho_B)-\rho_B}$ denotes the 
irreducible finite-dimensional
representation of $M$ with highest weight $\omega(\lambda_B+\rho_B)-\rho_B$.

\section{A Formula for $\Phi_M(\gm,\Theta^E)$}

To identify \eqref{Weyl} with  $\Phi_M(\gm,\Theta^E)$, we 
replace $m(\omega * 
B)$ with $n(\gm,\omega * B)$ as on page 500 of [2], and multiply it by the 
factor $\delta_P^{\half}(\gm)$:

\begin{equation}
\label{factor}
  \delta_P^{\half}(\gm) \cdot \sum_{\omega \in W^M} n(\gm,\omega * B) \cdot 
\varepsilon(\omega) \cdot 
\tr(\gm; V^M_{\omega(\lambda_B+\rho_B)-\rho_B}).
\end{equation}

\noindent Here $\delta_P$ is the modulus character of $P$.
(We are still only considering regular $\gm$.)

\bigskip
Write $A_G$ for the split component of the center of $G$.  Let
$\lambda_0 \in X^*(A_G)$ denote the character by which $A_G$ acts on $E$.  
It extends to $X^*(T)_{\R}$ in the usual way, and is 
$W$-invariant.
\bigskip

Let $T_e$ denote the subtorus of $T$ generated by $T_c$ and $A_G$.  It is 
the maximal subtorus of $T$ which is elliptic in $G$.
\bigskip

Write $p_M$ for the projection from $X^*(T)_{\R}$ to 
$X^*(A)_{\R}$, and note
that it is $W_L$-invariant.
The group $W_L$ fixes each root of $M$, thus it acts on $W^M$. 
For every orbit of this action, there is a unique 
member $\omega$ so that $p_M(\omega(\lambda_B + \rho_B-\lambda_0))$ is 
dominant with respect
to $C$.  We denote the set of these elements by $W^{LM}$, one element for 
each orbit of $W_L$ on $W_M$.

\bigskip
If $\lambda \in X^*(T)$ and $w_L \in W_L$, then 
plainly 
$w_L \lambda-\lambda \in \II_M^*$.  
Write $(\chi_{w_L,\omega,B},\CC_{w_L,\omega,B})$ for the 
one-dimensional 
representation of $M$, acting through $A$, with weight $w_L\omega(\lambda_B+\rho_B)-\omega(\lambda_B+\rho_B)$.
Note that $T_c$ and $A_G$ act trivially on 
$\CC_{w_L,\omega,B}$, thus so does $T_e$.

\bigskip
Thus we have
\[  V^M_{ w_L\omega(\lambda_B+\rho_B)-\rho_B} \cong V^M_{ 
\omega(\lambda_B+\rho_B)-\rho_B} \otimes \CC_{w_L,\omega,B}. \]

\bigskip
Our formula \eqref{factor} is now (replacing $\omega \in W^M$ 
with $w_L 
\omega$, where $\omega$ is now in $W^{LM}$):

\begin{equation}
\label{replace}
\delta_P^{\half}(\gm) \cdot \sum_{\omega \in W^{LM}} \varepsilon(\omega) \cdot 
\tr(\gm; V^M_{\omega(\lambda_B+\rho_B)-\rho_B}) \cdot
\sum_{w_L \in W_L}  \varepsilon(w_L) \cdot \chi_{w_L,\omega,B}(\gm) \cdot 
n(\gm, w_L \omega 
* B), 
\end{equation}

Of course we now wish to simplify the inner sum.  Recall from 
page 500 of [2] that 
\[ n(\gm,w_L \omega * B)= \bar{c}(x,p_M(w_L \omega \lambda_B+ w_L\omega 
\rho_B - \lambda_0)), 
\]

\noindent where $x$ is in the interior of $C$.  
Here $\bar{c}(x,\lambda)$ is the integer-valued ``stable discrete series 
constant'' on
\[ (X_*(A/A_G)_{\R})_{\reg} \times (X^*(A/A_G)_{\R})_{\reg}, \]

\noindent as defined, for instance, on page 493 of [2]. 
Recall that $\lambda_0 \in X^*(T)_{\R}$ is obtained from the character 
$\lambda_0 \in X^*(A_G)$ by which $A_G$ acts on $E$, and is thus 
$W$-invariant.

\bigskip

As $p_M$ commutes with $w_L$, the inner sum of 
\eqref{replace} is now

\begin{equation}
\label{inner}
 \sum_{w_L \in W_L}  \varepsilon(w_L) \cdot \bar{c}(x,w_L \Lambda) 
\cdot \chi_{w_L,\omega,B}(\gm), 
\end{equation}

\noindent where $\Lambda=p_M(\omega \lambda_B+ \omega \rho_B -
\lambda_0)$.  

\bigskip

We would like to consider the limit of \eqref{inner} as $x$ 
approaches 
$0$.   Recall we can write $\gm=\gm_c \cdot \exp(x)$, with 
$\gm_c \in T_c(\R)$ and $x$ in $\bar{C}$.  Also recall that $\gm$ is 
still regular (not for long!).  Consider the above formula with $\gm_c$ fixed and $x$ 
going to $0$ along regular elements of $\bar{C}$.  Fix some element $x_0$ in 
the interior of $C$.
The value 
\[ \bar{c}(x,w_L \Lambda)=\bar{c}(x_0,w_L \Lambda) \]
is unchanged, but $\chi_{w_L,\omega,B}(\gm)$ approaches 
$\chi_{w_L,\omega,B}(\gm_c)=1$.

Thus \eqref{inner} converges to 
\[ \sum_{w_L \in W_L}  \varepsilon(w_L) \cdot \bar{c}(x_0,w_L \Lambda) \]

\noindent for some $x_0 \in C$.

\bigskip

But this is simply $(-1)^{q(L)}|W_L|$, by  Proposition 1(ii) in 
Section 5 below.  Here we use that $\omega \in W^{LM}$.  Note that $-1$ is 
in the Weyl group of the root system by the argument on page 499 of [2].

\bigskip

It is easy to modify this argument to get the 
same limit as $x$ approaches an element of $X_*(A_G)_{\R}$.
\bigskip

Finally note that $\delta_P$ is a positive character and therefore trivial
on the compact group $T_c(\R)$.  It is thus trivial on $T_e(\R)$.

\bigskip

Now consider irregular $\gm$.  We take the limit in 
\eqref{replace} and obtain our theorem:

\begin{thm}

If $\gm \in T_e(\R)$, then 
\begin{equation}
\label{main}
\Phi_M(\gm,\Theta^E) =  (-1)^{q(L)}\cdot |W_L| \cdot \sum_{w \in 
W^{LM}} \varepsilon(w) \cdot \tr(\gm; V^M_{w(\lambda_B+\rho_B)-\rho_B}). 
\end{equation}
\end{thm}

\bigskip

\noindent For the reader's convenience, we review the definition of 
$W^{LM}$.

\bigskip

The definition depends on the choice of a parabolic $P=MN$ and a 
Borel subgroup $B$ with $T \subseteq B \subseteq P$.  The choice of $B$ gives a 
set of positive roots $R^+$ for $R$ and a set of positive roots $R_M^+$ 
for $R_M$.  It also gives $B$-dominant elements $\lambda_B$ and $\rho_B$
of $X^*(T)_{\R}$.  The choice of $P$ determines an $L$-chamber $C$ as in 
Section 2.  Recall the character $\lambda_0$ determined by $A_G$ on $E$ 
and the projection $p_M$ from $X^*(T)_{\R}$ to $X^*(A)_{\R}$.  Then 
\[ W^{LM} = \{ w \in W | w^{-1}R_M^+ \subseteq R^+ \text{ and } 
p_M(w(\lambda_B + \rho_B - \lambda_0)) \text{ is dominant w.r.t. } C \}. 
\]

\bigskip

We now evaluate \eqref{main} for $\Phi_M$ on the extreme cases 
for $T$.  If $T=A$ is split, then $M=A$, $L=G$, $W^{LM}$ is trivial, but so is
$T_c$.  We conclude that for $z \in A_G(\R)$,
\[ \Phi_A(z,\Theta^E)=(-1)^{q(G)} \cdot |W| \cdot \lambda_0(z). \]

If $T$ is elliptic, then $M=G$, $L=T$, $W^{LM}$ is again 
trivial, and so for $\gm \in T$,
\[ \Phi_G(\gm,\Theta^E)= \tr(\gm; E). \]

\noindent Note that this agrees with the results of 
Theorems 5.1 and 5.2 of [2], since
\[ \tr(\gm^{-1}; E^*)= \tr(\gm; E). \]

\section{The Sum of the Stable Discrete Series Constants}

Let $(X,X^*,R,\check{R})$ be a root system.  Write $W$ for the Weyl group of the root 
system, and $\varepsilon$ for its sign character.  Assume that $R$ 
generates the real vector space $X$ and that $-1 \in W$.  
Write $q(R)$ for $(|R^+|+\dim(X)) / 2$, as in [2].
Let $x_0$ be a regular element of $X$, and $\lambda$ a regular 
element of $X^*$.  Write $C_0$ for the chamber of $X$ containing 
$x_0$, and $C_0^{\vee}$ for its dual chamber in $X^*$.
Recall the stable discrete series constants $\bar{c}_R(x_0,\lambda)$ from 
section 3 of [2]. 

\begin{prop} We have the following formulas for sums of discrete 
series constants:  
\begin{itemize}

\item[(i)] For all such  $ \lambda, \sum_{w \in W} 
\bar{c}_R(wx_0,\lambda)=|W|.$ 

\item[(ii)] For $\lambda=\lambda_0 \in C_0^{\vee}$, we have $\sum_{w 
\in W} 
\varepsilon(w) \cdot \bar{c}_R(wx_0,\lambda_0)=(-1)^{q(R)}|W|.$

\end{itemize}

\noindent The same formulas hold if we sum over the $W$-orbit of 
$\lambda$ rather than
that of $x_0$.

\end{prop}

We make a few comments before beginning the proof.  The proof 
begins by using the   
``inductive'' property (4) of the discrete series constants from page 493 of 
[2], to change the sum 
over chambers into a sum over certain facets of $X$.  In fact we consider those facets which  
separate the chambers of $X$, i.e., those which span the root hyperplanes $Y$ of $X$.

\bigskip

In the course of the proof, we (mis)use the term ``facet'' only 
in reference to these 
particular facets, of codimension $1$.  So a facet in this sense will be the common face 
of two adjacent chambers.

\bigskip

The hyperplanes $Y$ have their own chambers, and we examine the 
relationship between the 
facets and these smaller chambers.  Not every facet is equal to such a chamber, as in the 
case of $B_3$ when $Y$ is the root hyperplane of a long root.  The facets 
in $Y$ give a $B_2$ system, but the chambers of $R_Y$ give an $A_1 
\times A_1$ system.

\bigskip

Finally induction on the rank of the root system gives the 
calculation.

\begin{proof}
\noindent The second formula follows from the first by applying Theorem 
3.2(2) on page 494 of [2].

\bigskip

We induce on $r=\dim X$.  The proposition is clear when $r=0$.

\bigskip

We associate these discrete series constants with the 
various chambers and 
facets of $X$, and introduce some appropriate notation.

\bigskip

Write $c(\Ch)$ for $\bar{c}_R(x,\lambda)$, when $x$ is in the 
interior of a chamber 
$\Ch$.

\bigskip

Suppose $F$ is a facet in $X$, $y$ is in the interior of $F$, 
and 
$\bar{F}:=\Span(F)=Y$.  
Then write $c(F)=\bar{c}_{R_Y}(y,\lambda_Y)$, notation as on page 493 of 
[2].

\bigskip

Thus if $F$ is the common face of distinct chambers $\Ch$ and 
$\Ch^{\prime}$, then
\[ 2c(F)=c(\Ch)+c(\Ch^{\prime}). \]

Each chamber has $r$ faces, and it follows that 

\begin{equation}
\label{summing}
  r \cdot \sum_{\Ch} c(\Ch)= 2 \sum_{F} c(F), 
\end{equation}
\noindent where we are summing over all chambers and then all facets.

\bigskip

We show the right hand side of \eqref{summing} is equal to $r 
\cdot |W|$ to 
prove the proposition.

\bigskip

Now every facet is on some root hyperplane 
$X_{\alpha}=X_{-\alpha}$, so we have
\[ 2 \sum_F c(F)= \sum_{\alpha \in R} \sum_{\bar{F}=X_{\alpha}} c(F). \]

We now work with the inner sum.  There is a root system on 
$X_{\alpha}$ whose set of coroots 
is $\check{R} \cap X_{\alpha}$, which defines chambers $\Ch_{\alpha}$ in $X_{\alpha}$ and 
constants $c_{\alpha}(\Ch_{\alpha})$.  Write $W_{\alpha}$ for 
the Weyl group of $X_{\alpha}$.  We have
\[ \sum_{\bar{F}=X_{\alpha}} c(F)=\sum_{\Ch_{\alpha}} \sum_{F \subset 
\Ch_{\alpha}} c_{\alpha}(\Ch_{\alpha})= 
\sum_{\Ch_{\alpha}} \sum_{W_{\alpha} \backslash \{F \subset 
X_{\alpha} \}} c_{\alpha}(\Ch_{\alpha}) = \sum_{W_{\alpha} \backslash \{ F \subset
X_{\alpha} \} } \sum_{\Ch_{\alpha}}c_{\alpha}(\Ch_{\alpha}). \]

\noindent For the first equality, note that every facet $F$ with 
$\bar{F}=X_{\alpha}$ is contained in 
a some chamber $\Ch_{\alpha}$. 

\bigskip

The second equality follows because $W_{\alpha}$ acts transitively on the 
chambers $C_{\alpha}$.

\bigskip

Write $n(\alpha)$ for the order of $W_{\alpha} \backslash \{ F 
\subset X_{\alpha} \}$.  
It is equal to the number of facets in a given chamber $\Ch_{\alpha}$.  Then 
by induction the above is merely
\[ n(\alpha) \cdot |W_{\alpha}|, \]

\noindent which is exactly the number of facets in $X_{\alpha}$.  It 
follows that \eqref{summing} is 
simply equal to twice the total number of facets in $X$.

\bigskip

Since $W$ has $r$ orbits on the set of facets in $X$, and the 
stabilizer in $W$ of any 
facet has order $2$, we conclude that the total number of facets is half of $r \cdot |W|$, 
as desired.
\end{proof}

\end{document}